\documentclass[a4paper,12pt]{article}
\usepackage{amssymb,amsthm,amsmath,fullpage,graphicx}

\newtheorem{thm}{Theorem}[section]
\newtheorem{defn}[thm]{Definition}

\newtheorem{cor}[thm]{Corollary}

\newtheorem{lem}[thm]{Lemma}
\newtheorem{propn}[thm]{Proposition}
\newtheorem{rem}[thm]{Remark}

\begin{document}

\title{Moduli of continuity of local times\\of random walks on graphs\\in terms of the resistance metric}
\author{D.~A.~Croydon}
\maketitle

\begin{abstract}
In this article, universal concentration estimates are established for the local times of random walks on weighted graphs in terms of the resistance metric. As a particular application of these, a modulus of continuity for local times is provided in the case when the graphs in question satisfy a certain volume growth condition with respect to the resistance metric. Moreover, it is explained how these results can be applied to self-similar fractals, for which they are shown to be useful for deriving scaling limits for local times and asymptotic bounds for the cover time distribution.
\end{abstract}

\section{Introduction}

Over the last couple of decades, extensive efforts have been devoted to studying the behaviour of random walks on general graphs, work that has yielded, for instance, estimates for the corresponding heat kernel and mixing times in terms of quantities such as volume growth and electrical resistance, which do not depend on precise structural information (see, for example, \cite{BCK, BJKS, KM, NP}). Furthermore, for a wide range of families of fractal and random graphs, scaling limits have been established for the laws of the random walks upon them \cite{BP, CCK, Croy1, Croy3, Croy4, K, KK}, as well as corresponding asymptotic results for heat kernels and mixing times \cite{CH,CHK}. More delicate properties of random walks on graphs, particularly the cover time, are also becoming better understood. Indeed, recent years have seen the order of growth of the cover time computed for various families of graphs \cite{Abe, BDNP}, and a strong connection has been made between the cover time and the maximum of the Gaussian free field for any graph \cite{DLP}. Moreover, in some special cases where there is concentration of the cover time about its mean, extremely precise distributional convergence results are known, notably for the two-dimensional discrete torus \cite{DPRZ, Ding}. Partly motivated by providing techniques for studying the cover time in settings where there is not concentration of the cover time -- as is the case for many self-similar fractals, in this article we study the continuity properties of local times on graphs. Since the first time that local times of a simple random walk on a graph are non-zero everywhere gives the cover time, we believe that our results will provide another tool for studying the latter; this is a point upon which we will expand later in the article.

In order to present our main results, let us start by introducing the framework in which we are working. In particular, let $G=(V(G),E(G))$ be a finite connected graph, where $V(G)$ denotes the vertex set and $E(G)$ the edge set of $G$. To avoid trivialities, we always assume that $G$ has at least two vertices. Let $\mu^G:V(G)^2\rightarrow\mathbb{R}_+$ be a weight function that is symmetric, i.e.~$\mu_{xy}^G=\mu_{yx}^G$, and satisfies $\mu_{xy}^G>0$ if and only if $\{x,y\}\in E(G)$. The associated discrete time simple random walk is then the Markov chain $((X^G_t)_{t\geq 0}, \mathbf{P}^G_x,x\in V(G))$ with transition probabilities $(P_G(x,y))_{x,y\in V(G)}$ defined by
\[P_G(x,y):=\frac{\mu^G_{xy}}{\mu^G_{x}},\]
where $\mu^G_x:=\sum_{y\in V(G)}\mu^G_{xy}$. We note that the invariant probability measure of this process is a multiple of the measure version of $\mu^G$ obtained by setting $\mu^G(\{x\}):=\mu_x^G$ for $x\in V(G)$. The process $X^G$ has corresponding local times $(L^G_t(x))_{x\in V(G),t\geq 0}$, given by $L^G_0(x)= 0$ and, for $t\geq 1$,
\[L^G_t(x):=\frac{1}{\mu_x^G}\sum_{i=0}^{t-1}\mathbf{1}_{\{X^G_i=x\}}.\]
It is providing a modulus of continuity of these random functions in the spatial variable $x$ that is the focus of this article.

It is widely known that there are close connections between the study of random walks on graphs and electrical networks. Our work will contribute to this area by providing estimates for the fluctuations in the local times of a random walk in terms of the so-called resistance metric, which we now introduce. More specifically, the process $X^G$ has an associated Dirichlet form given by
\[\mathcal{E}_G(f,g):=\frac{1}{2}\sum_{\substack{x,y\in V(G):\\ \{x,y\}\in E(G)}}(f(x)-f(y))(g(x)-g(y))\mu_{xy}^G,\]
for $f,g:V(G)\rightarrow \mathbb{R}$, which can in turn be used to define the resistance operator through the variational formula
\[R_G(A,B)^{-1}:=\inf\left\{\mathcal{E}_G(f,f):\:f:V(G)\rightarrow\mathbb{R},\:f|_A=0,\:f|_B=1\right\}\]
for $A,B$ disjoint subsets of $V(G)$; the latter is so-called because it describes the effective resistance between $A$ and $B$ in the graph when it is viewed as an electrical network with conductances along edges given by the weight function $\mu^G$. We then define the resistance metric on the vertices of $G$ by setting $R_G(x,y):=R_G(\{x\},\{y\})$ if $x\neq y$, and $R_G(x,x):=0$. We note that the resistance metric is indeed a metric (see \cite[Proposition 4.25]{Barlow}).

In studying cover times of random walks on graphs, continuity properties of local times in terms of the resistance metric have previously been considered. Indeed, applying the key identity
\begin{equation}\label{returnprob}
\mathbf{P}^G_x\left(\tau_y<\tau_x^+\right)=\frac{1}{\mu_x^GR_G(x,y)},
\end{equation}
where $\tau_y=\inf\{t\geq 0:\:X^G_t=y\}$ is the first hitting time of $y$, and $\tau_x^+=\inf\{t\geq 1:\:X^G_t=x\}$ the first return time to $x$ (see \cite[Proposition 9.5]{LPW}, for instance), the following Gaussian concentration result was established as \cite[Lemma 5.2]{KKLV}: if $\tau_x(0)=0$ and, for $i\geq 1$, $\tau_x(i)$ is defined to be the time of the $i$th subsequent visit to $x$ by $X^G$ (so that $\tau_x$ can be considered to be the inverse local time at $x$), then
\[\mathbf{P}^G_x\left(L^G_{\tau_x(i)}(x)-L^G_{\tau_x(i)}(y)\geq \lambda\right)\leq e^{- \frac{\lambda^2\mu_x^G}{4iR_G(x,y)}},\]
for all $i\geq 0$, $\lambda>0$ and $x,y\in V(G)$ (cf. \cite[Lemma 1.12]{DLP}). Our concentration estimates (Theorem \ref{mainest} below) are of a similar form, with an important distinction being that we are interested in estimating the fluctuations of the local time at deterministic times, rather than at the random inverse local time. In the first estimate, we sacrifice Gaussian tails to establish a bound that holds uniformly over time intervals -- this is the discrete time analogue of \cite[(V.3.28)]{BG}. In the second, Gaussian tails are obtained at the cost of truncating the local times (cf.~\cite[Lemma 2.8]{Barlowcont}). For the statement of these bounds, as two important measures of the scale of $G$, we define
\[m(G):=\mu^G(V(G)),\qquad r(G):=\max_{x,y\in V(G)}R_G(x,y),\]
to be its total mass with respect to the measure $\mu^G$, and its diameter in the resistance metric, respectively. Note that the product $m(G)r(G)$ gives the maximal commute time of the random walk, e.g.~\cite[Proposition 10.6]{LPW}, and so gives a natural time-scaling. We also introduce the rescaled resistance metric $\tilde{R}_G(x,y):=r(G)^{-1}{R}_G(x,y)$, and define the notation $x\wedge y:=\min\{x,y\}$.

{\thm \label{mainest} (a) For each $T>0$, there exist constants $c_1$ and $c_2$ not depending on $G$ such that
\[\max_{x,y,z\in V(G)}\mathbf{P}^G_z\left(\max_{0\leq t\leq Tm(G)r(G)}r(G)^{-1}\left|L^G_t(x)-L^G_t(y)\right|\geq \lambda\sqrt {\tilde{R}_G(x,y)}\right)\leq c_1e^{-c_2 \lambda}\]
for every $\lambda\geq 0$. (NB. The constants can be chosen such that only $c_1$ depends on $T$.)\\
(b) For any $G$ it holds that
\[\max_{x,y,z\in V(G)}\mathbf{P}^G_z\left(\max_{t\geq 0}\left|L\wedge \left(\frac{L^G_t(x)}{r(G)}\right)-L\wedge \left(\frac{L^G_t(y)}{r(G)}\right)\right|\geq \lambda\sqrt {\tilde{R}_G(x,y)}\right)\leq 2e^{\frac{1}{2}-\frac{\lambda^2}{8L}}\]
for every $\lambda\geq 0$ and $L\geq 1$.}
\bigskip

To provide a modulus of continuity for the local times, the goal is to bring the maximum over the arguments of the local times inside the estimates of the previous result. In the continuous setting, this has been achieved by applying a general estimate on the fluctuations of a function on Euclidean space known in the literature as Garsia's lemma (after \cite[Lemma 1]{Garsia}, cf.~\cite[Lemma 1.1]{GRR}). Indeed, this approach was first used in \cite[Theorem 2]{GK} to deduce the continuity of local times of Markov processes on the real line. (A similar argument was applied to the local times of the Brownian motion on the Sierpi\'nski gasket in \cite[Theorem 1.11]{BP}.) Moreover, the argument was subsequently strengthened for L\'{e}vy processes in \cite{Barlowcont} (see also the estimates for the Sierpi\'nski carpet that appear as \cite[Theorem 8.2]{BB}). The aim here is to adapt the same approach to the discrete setting, and so for this we derive below a version of Garsia's lemma for graphs (see Proposition \ref{garsia}). To obtain a modulus of continuity estimate from this, the one restriction we need is some uniform control on the volume growth of the graphs in question. To state the condition we need, we define
\[B_G(x,r):=\left\{y\in V(G):\:R_G(x,y)<r\right\}\]
to be the open ball in the resistance metric, and set $r_0(G):=\min_{x,y\in V(G):x\neq y}R_{G}(x,y)$.

{\defn A collection of finite connected weighted graphs $(G_i)_{i\in I}$ is said to satisfy \emph{uniform volume growth with volume doubling (UVD)} if there exist constants $c_1,c_2,c_3\in(0,\infty)$ such that
\[c_1v(r)\leq \mu^{G_i}\left(B_{G_i}(x,r)\right)\]
for every $x\in G_i$, $r\in [r_0(G_i),r(G_i)]$, $i\in I$, and moreover,
\[m(G_i)\leq c_2v(r(G_i))\]
for every $i\in I$, where $v:\mathbb{R}_+\rightarrow\mathbb{R}_+$ is non-decreasing function with $v(2r)\leq c_3v(r)$ for every $r\in\mathbb{R}_+$.}

{\rem The above condition is weaker than is sometimes called uniform volume growth with volume doubling, since we do not require the upper volume bound to hold for balls smaller than the full space.}
\bigskip

Our next main result is that, under UVD, $\tilde{R}_{G_i}(x,y)^{1/2}(1+\ln \tilde{R}_{G_i}(x,y)^{-1})^{1/2}$ provides, with uniformly high probability, a modulus of continuity for the rescaled local times $r(G_i)^{-1}L^{G_i}_t(x)$ in the spatial variable (uniformly over the appropriate time interval). Observe that the particular form of $v$ that appears in the UVD property does not affect the modulus of continuity.

{\thm \label{mainresult} If $(G_i)_{i\in I}$ is a collection of graphs that satisfies UVD, then, for each $T>0$,
\[\lim_{\lambda\rightarrow\infty}\sup_{i\in I}\max_{z\in V(G_i)}\mathbf{P}^{G_i}_z\left(\max_{x,y\in V(G_i)}\max_{0\leq t\leq Tm(G_i)r(G_i)}\frac{r(G_i)^{-1}\left|L^{G_i}_{t}(x)-L_t^{G_i}(y)\right|}{\sqrt{\tilde{R}_{G_i}(x,y)\left(1+\ln \tilde{R}_{G_i}(x,y)^{-1}\right)}}\geq \lambda\right)=0.\]}
\bigskip

Whilst we do not pursue it here, we expect that the above bound on the modulus of continuity for rescaled local times is sharp up to constants. In the one-dimensional case, where $X^{G_i}$ is simple random walk on the interval $\{0,1,\dots,i\}$, in particular, it should be possible to check this, either directly or by a coupling of the rescaled random walks with reflected Brownian motion on the interval, using a Ray-Knight-type description of the relevant local times and an appeal to Levy's modulus of continuity theorem for Brownian motion (or the discrete adaptation thereof).

The remainder of the article is organised as follows. In Section \ref{ltsec}, we establish the concentration estimates of Theorem \ref{mainest}. Our discrete version of Garsia's lemma is established in Section \ref{garsec}, and applied under the assumption of UVD in Section \ref{uvdsec}, thereby proving Theorem \ref{mainresult}. Subsequently, in Section \ref{examsec}, we present a number of examples, including self-similar fractal graphs, to which these results apply. Moreover, in Section \ref{infsec}, an adaptation of the results to a class of infinite graphs is derived. Finally, in Section \ref{scalingsec}, we consider some of the consequences of our equicontinuity results. In particular, we show that if we have a sequence of graphs such that the associated random walks admit a diffusion scaling limit that has jointly continuous local times, and a suitable local time equicontinuity result holds, then it is further possible to obtain convergence of rescaled local times. We also discuss an application to the study of cover times of random walks on graphs.

\section{Local time concentration estimates}\label{ltsec}

The aim of this section is to prove Theorem \ref{mainest}. We begin with a lemma that provides an estimates for the distributional tail of the return time. We note that similar estimates to this and the next result have appeared elsewhere in the literature, e.g.~\cite[Section 3]{KP}, but we include their proofs for completeness.

{\lem \label{returntail} There exist universal constants $c_1,c_2$ such that
\[\mu_x^Gr(G)\mathbf{P}_x^G\left(\tau_x^+\geq \lambda m(G)r(G)\right)\leq c_1e^{-c_2\lambda}\]
for every $\lambda\geq 4$.}
\begin{proof} By applying the Markov property at times $2\lceil m(G)r(G)\rceil ,4\lceil m(G)r(G)\rceil ,\dots$, we deduce that $\mathbf{P}_x^G(\tau_x^+\geq \lambda m(G)r(G))$ is bounded above by
\[\mathbf{P}_x^G\left(\tau_x^+\geq 2\lceil m(G)r(G)\rceil\right)\left(\max_{y\in  V(G)\backslash\{x\}}\mathbf{P}_y^G\left(\tau_x\geq 2\lceil m(G)r(G)\rceil\right)\right)^{k-1},\]
where $k:=\lfloor \lambda m(G)r(G)/2\lceil m(G)r(G)\rceil\rfloor\geq \lfloor \lambda/4\rfloor \geq 1$ (note that $m(G)r(G)\geq 1$). Since
\begin{equation}\label{mom}
\mathbf{E}_x^G\tau_x^+=\frac{m(G)}{\mu_x^G},\qquad \mathbf{E}_y^G\tau_x\leq R_G(x,y)m(G)
\end{equation}
for all $x,y\in V(G)$ (see, for example, \cite[Lemma 10.5 and Proposition 10.6]{LPW}), we obtain that
\[\mathbf{P}_x^G\left(\tau_x^+\geq \lambda m(G)r(G)\right)\leq
\frac{m(G)}{2\mu_x^G\lceil m(G)r(G)\rceil}\left(\frac{r(G)m(G)}{2\lceil m(G)r(G)\rceil}\right)^{k-1}\leq  \frac{1}{\mu_x^G r(G)}\left(\frac12\right)^{k},\]
and the result follows.
\end{proof}

The above lemma readily yields the following corollary.

{\cor \label{mgfcor}There exists a universal constant $c$ such that
\[\mathbf{E}_x^G\left(e^{-\theta\tau_x^+}\right)\leq e^{- \frac{\theta m(G)}{\mu_x^G}+ \frac{c\theta^2 m(G)^2r(G)}{\mu_x^G}}\]
for every $\theta\geq 0$.}
\begin{proof} Since $1-x\leq e^{-x}\leq 1-x+x^2/2$ for $x\geq 0$, we have that
\begin{equation}\label{m1}
\mathbf{E}_x^G\left(e^{-\theta\tau_x^+}\right)\leq 1- \theta\mathbf{E}_x^G\left(\tau_x^+\right)+\frac{\theta^2}2\mathbf{E}_x^G\left((\tau_x^+)^2\right)\leq e^{-\theta\mathbf{E}_x^G\left(\tau_x^+\right)+\frac{\theta^2}2\mathbf{E}_x^G\left((\tau_x^+)^2\right)}.
\end{equation}
From (\ref{mom}), we know that $\mathbf{E}_x^G\left(\tau_x^+\right)=m(G)/\mu_x^G$. For the second moment, we apply this and Lemma \ref{returntail} to deduce
\begin{eqnarray*}
\mathbf{E}_x^G\left((\tau_x^+)^2\right)&\leq &
2\sum_{k=0}^\infty k \mathbf{P}_x^G\left(\tau_x^+\geq k\right)\\
&\leq & 2\sum_{k=1}^{4m(G)r(G)}\mathbf{E}_x^G\left(\tau_x^+\right)+\frac{2c_1}{\mu_x^Gr(G)}\sum_{k=4m(G)r(G)+1}^\infty ke^{-c_2 k/m(G)r(G)}\\
&\leq &\frac{8m(G)^2r(G)}{\mu_x^G}+\frac{2c_1e^{-c_2 /m(G)r(G)}}{\mu_x^Gr(G)\left(1-e^{-c_2/m(G)r(G)}\right)^2}\\
&\leq & \frac{cm(G)^2r(G)}{\mu_x^G}.
\end{eqnarray*}
For the final inequality, we again use that $m(G)r(G)\geq1$. Inserting this estimate into (\ref{m1}), we obtain the desired result.
\end{proof}

\begin{proof}[Proof of Theorem \ref{mainest}(a).] To begin with, suppose that the random walk starts from $z=x$, where $x\neq y$. As in the introduction, let $\tau_x(0)=0$ and, for $i\geq 1$, $\tau_x(i)$ be the time of the $i$th subsequent visit to $x$ by $X^G$. Since local times are monotonic in the time variable, we have for $t\in (\tau_x(i),\tau_x(i+1)]$ that $|L^G_t(x)-L^G_t(y)\leq |L^G_t(x)-L^G_{\tau_x(i)}(y)|+|L^G_t(x)-L^G_{\tau_x(i+1)}(y)|$. For $t$ in this range, we also have that $L^G_t(x)=L^G_{\tau_x(i+1)}(x)$. Moreover, $L^G_{\tau_x(i+1)}(x)-L^G_{\tau_x(i)}(x)=1/\mu_x^G$. Thus
\[\left|L^G_t(x)-L^G_t(y)\right|\leq \sum_{j=0}^1\left|L^G_{\tau_x(i+j)}(x)-L^G_{\tau_x(i+j)}(y)\right|+\frac{1}{\mu_x^G}.\]
Consequently, for $\lambda>0$ and $L\geq 0$,
\begin{eqnarray*}
\lefteqn{\mathbf{P}^G_x\left(\max_{0\leq t\leq Tm(G)r(G)}\left|L^G_t(x)-L^G_t(y)\right|\geq \lambda \sqrt {r(G)R_G(x,y)}\right)}\\
&\leq &\mathbf{P}^G_x\left(\mu^G_x L^G_{ Tm(G)r(G)}(x)> L+1\right)+\mathbf{P}^G_x\left(1/\mu_x^G\geq \lambda \sqrt { r(G)R_G(x,y)}/2\right)\\
&&+\mathbf{P}^G_x\left(\max_{1\leq i\leq {L+1}}
\left|L^G_{\tau_x(i)}(x)-L^G_{\tau_x(i)}(y)\right|\geq \lambda \sqrt { r(G)R_G(x,y)}/4\right)\\
&=:&T_1+T_2+T_3,
\end{eqnarray*}
where we note that the condition $\mu^G_x L^G_{ Tm(G)r(G)}(x)>L+1$ is equivalent to $\tau_x(L+1)< Tm(G)r(G)$. We will bound each of the above three terms separately. To this end, we first note that $\mu_x^G\sqrt{ r(G)R_G(x,y)}\geq 1$, and so $T_2$ is equal to 0 whenever $\lambda > 2$.

We next consider $T_1$. Using the fact that under $\mathbf{P}_x^G$ the variables $(\tau_{x}(i+1)-\tau_{x}(i))_{i\geq 0}$ form an independent sequence, each distributed as $\tau_x^+$, one can deduce that, for $\theta\geq 0$,
\begin{eqnarray*}
T_1&=& \mathbf{P}_x^G\left(\tau_{x}(L+1)< { Tm(G)r(G)}\right)\\
&=&\mathbf{P}_x^G\left(\sum_{i=0}^{L}\left(\tau_{x}(i+1)-\tau_{x}(i)-\frac{m(G)}{\mu_x^G}\right)< { Tm(G)r(G)}-\frac{(L+1)m(G)}{\mu_x^G}\right)\\
&\leq & e^{-\theta \left(\frac{(L+1)m(G)}{\mu_x^G}-Tm(G)r(G)\right)}\mathbf{E}_x^G\left(e^{-\theta\left(\tau_x^+-\frac{m(G)}{\mu_x^G}\right)}\right)^{L+1}.
\end{eqnarray*}
By Corollary \ref{mgfcor}, it therefore holds that
\[T_1 \leq e^{-\theta \left(\frac{(L+1)m(G)}{\mu_x^G}-Tm(G)r(G)\right)+ \frac{c(L+1)\theta^2 m(G)^2r(G)}{\mu_x^G}},\]
and optimising over $\theta$ yields, for $L+1\geq T\mu_x^Gr(G)$,
\[T_1 \leq e^{-\mu_x^G\left(\frac{(L+1)}{\mu_x^G}-Tr(G)\right)^2/2c(L+1)r(G)}.\]

We now turn to $T_3$, and for the moment we assume that $\mu_y^GR_G(x,y)>1$. Observe that the term in question can be written as
\[T_3=\mathbf{P}^G_x\left(\max_{1\leq i\leq L+1}\left|\sum_{j=1}^{i} \left(\frac{1}{\mu_x^G}-\eta_j\right)\right|\geq \lambda \sqrt { r(G) R_G(x,y)}/2\right),\]
where $\eta_i:=L^G_{\tau_x(i)}(y)-L^G_{\tau_x(i-1)}(y)$ for $i\geq 1$.
Now, $(\eta_i)_{i\geq 1}$ is an independent and identically distributed sequence, and it is a simple application of the Markov property to obtain from (\ref{returnprob}) that if $N_i:=\mu^G_y\eta_i$, then
\begin{equation}\label{nprob}
\mathbf{P}^G_x\left(N_i=k\right)=\frac{1}{\mu^G_xR_G(x,y)}\left(1-\frac{1}{\mu_y^G R_G(x,y)}\right)^{k-1}\frac{1}{
\mu_y^GR_G(x,y)},
\end{equation}
for $k\geq 1$. Moreover, it is elementary to check from (\ref{nprob}) that $\mathbf{E}^G_x(\eta_i)=(\mu_x^G)^{-1}$. This means that $(M_i)_{i\geq 0}$, where $M_0=0$ and, for $i\geq1$,
\begin{equation}\label{mdef}
M_i:=\sum_{j=1}^{i} \left(\frac{1}{\mu_x^G}-\eta_i\right),
\end{equation}
is a martingale, and further, for
$\theta\in(0,-\mu_y^G\ln(1-1/\mu_y^GR_G(x,y)))$, we have that $(e^{\theta |M_i|})_{i\geq 0}$ is a sub-martingale. (The condition $\theta <-\mu_y^G\ln(1-1/\mu_y^GR_G(x,y))$ ensures integrability.) Therefore, applying Doob's sub-martingale inequality, we deduce that
\begin{eqnarray*}
T_3&=&\mathbf{P}^G_x\left(\max_{0\leq i\leq L+1}\left|M_i\right|\geq \lambda \sqrt { r(G)R_G(x,y)}/2\right)\\
&\leq& \mathbf{P}^G_x\left(\max_{0\leq i\leq L+1}e^{\theta |M_i|}\geq e^{\theta\lambda \sqrt {r(G)R_G(x,y)}/2}\right)\\
&\leq &\mathbf{E}^G_x\left(e^{\theta |M_{L+1}|}\right)e^{-\theta\lambda \sqrt { r(G)R_G(x,y)}/2}\\
&\leq & \left(\mathbf{E}^G_x\left(e^{\theta M_{L+1}}\right)+\mathbf{E}^G_x\left(e^{-\theta M_{L+1}}\right)\right)e^{-\theta\lambda \sqrt {r(G)R_G(x,y)}/2}\\
&=&\left(\mathbf{E}^G_x\left(e^{-\theta\left(\eta_i-1/{\mu_x^G}\right)}\right)^{L+1}+
\mathbf{E}^G_x\left(e^{\theta\left(\eta_i-1/{\mu_x^G}\right)}\right)^{L+1}\right)e^{-\theta\lambda \sqrt { r(G)R_G(x,y)}/2}.
\end{eqnarray*}
A routine computation using (\ref{nprob}) gives that $\mathbf{E}^G_x(e^{\theta\left(\eta_i-1/{\mu_x^G}\right)})$ is equal to \[e^{-\theta/\mu_x^G}\left(1-\frac{1}{\mu_x^GR_G(x,y)}\right)
+\frac{e^{\theta(1/\mu_y^G-1/\mu_x^G)}}{\mu_x^GR_G(x,y)^2\mu_y^G\left(1-e^{\theta/\mu_y^G}\left(1-\frac{1}{\mu_y^GR_G(x,y)}\right)\right)}.\]
By considering the Taylor expansion of this expression, we deduce that
\[\mathbf{E}^G_x(e^{\theta\left(\eta_i-1/{\mu_x^G}\right)})\leq 1+\theta^2
\mathbf{E}^G_x\left(\left(\eta_i-1/{\mu_x^G}\right)^2\right)\leq e^{\theta^2
\mathbf{E}^G_x\left(\left(\eta_i-1/{\mu_x^G}\right)^2\right)}\]
uniformly over $\theta\leq c_1\min\{\mu_x^G,\mu_y^G,R_G(x,y)^{-1}\}=c_1R_G(x,y)^{-1}$, where $c_1\in(0,1)$ is some small universal constant. (Note that $R_G(x,y)^{-1}\leq -\mu_y^G\ln(1-1/\mu_y^GR_G(x,y))$, and so the integrability condition for the martingale is also satisfied if $\theta\in(0,c_1R_G(x,y)^{-1}]$.) Again appealing to (\ref{nprob}), it is possible to compute that
\begin{equation}\label{2ndmom}
\mathbf{E}^G_x\left(\left(\eta_i-1/{\mu_x^G}\right)^2\right)=\frac{2\left(1-\frac{1}{\mu_y^GR_G(x,y)}\right)R_G(x,y)}{\mu_x^G}+\frac{1}{\mu_x^G\mu_y^G}-\frac{1}{(\mu_x^G)^2}\leq \frac{2R_G(x,y)}{\mu_x^G}.
\end{equation}
So we obtain
\[T_3\leq 2e^{2\theta^2(L+1)R_G(x,y)/\mu_x^G -\theta\lambda\sqrt {r(G) R_G(x,y)}/2 }.\]
Again optimising over $\theta$, we find
\[T_3\leq  2e^{-\mu_x^G\lambda^2r(G) /16(L+1)},\]
at least assuming that $\mu_x^G\lambda\sqrt{(r(G)R_G(x,y))}/8(L+1)\leq c_1$.

In summary, we have so far shown that if $\mu_y^GR_G(x,y)>1$, $\lambda>2$ and it also holds that $L+1\geq T\mu_x^Gr(G)$ and $\mu_x^G\lambda\sqrt{(r(G)R_G(x,y))}/8(L+1)\leq c_1$, then
\begin{eqnarray*}
\lefteqn{\mathbf{P}_x^G\left(\max_{0\leq t\leq Tm(G)r(G)}\left|L^G_t(x)-L^G_t(y)\right|\geq \lambda \sqrt{r(G)R_G(x,y)}\right)}\\
&\leq&
e^{-\mu_x^G\left(\frac{(L+1)}{\mu_x^G}-Tr(G)\right)^2/2c(L+1)r(G)}+
 2e^{-\mu_x^G\lambda^2r(G)/16(L+1)}.
\end{eqnarray*}
Setting $L+1:=c_1^{-1}\lambda \mu_x^Gr(G)$, so that
\[\frac{\mu_x^G\lambda\sqrt{{r(G)}{R_G(x,y)}}}{8(L+1)}=\frac{c_1}{8}\sqrt{\frac{R_G(x,y)}{r(G)}}\leq c_1\]
we obtain that the relevant probability is bounded above by  $c_2e^{-c_3\lambda}$ for $\lambda>\max\{2,Tc_1\}$, where only $c_2$ depends on $T$. This bound is readily extended to hold for all $\lambda\geq 0$ by adjusting the constants as necessary, which establishes the desired estimate in this case.

When $\mu_y^GR_G(x,y)=1$, essentially the same argument applies. The main difference is that, because in this case the vertex $y$ is connected to $x$ by a single edge of resistance $R_G(x,y)=1/\mu_y^G$, the distribution of $N_i$ is given by
\[\mathbf{P}_x^G(N_i=0)=1-\frac{1}{\mu_x^GR_G(x,y)}=1-\mathbf{P}_x^G(N_i=1).\]
In particular, the bound for $T_1$ does not change, $T_2$ is still equal to 0 for $\lambda>2$, and a similar martingale argument can be used to estimate $T_3$. We omit the details.

Clearly, one could also reverse the role of $x$ and $y$ in the above argument, so that we start the process $X^G$ from $y$ instead. Hence, in the case that we start from an arbitrary vertex $z$, by applying the strong Markov property at the first time we hit the set $\{x,y\}$ (noting that the local times of both $x$ and $y$ are zero up until this time), one can also deduce a similar result. This concludes the proof.
\end{proof}

\begin{proof}[Proof of Theorem \ref{mainest}(b).] We note that the proof of this part of the theorem is an adaptation of the proof of \cite[Lemma 2.8]{Barlowcont}. We will use the same notation as in the proof of Theorem \ref{mainest}(a), though to account for arbitrary starting points will redefine $\tau_x(0)$ to be the first hitting time of $x$. (For $i\geq 1$, $\tau_x(i)$ will continue to denote the time of the $i$th subsequent visit to $x$ by $X^G$). In particular, for $x,y,z$ and $x\neq y$, by applying the strong Markov property at $\tau_x(0)$, this allows us to deduce that, for $L\geq 0$,
\begin{eqnarray*}
\lefteqn{\mathbf{P}^G_z\left(\max_{0\leq t\leq \tau_x(L+1)}\left(L^G_t(x)-L^G_t(y)\right)\geq \lambda \sqrt {r(G)R_G(x,y)}\right)}\\
&\leq &\mathbf{P}^G_x\left(\max_{0\leq t\leq \tau_x(L+1)}\left(L^G_t(x)-L^G_t(y)\right)\geq \lambda \sqrt {r(G)R_G(x,y)}\right)\\
&=&\mathbf{P}^G_x\left(\max_{1\leq i\leq L+1}\left(\frac{i+1}{\mu_x^G}-L^G_{\tau_x(i)}(y)\right)\geq \lambda \sqrt {r(G)R_G(x,y)}\right)\\
&=&\mathbf{P}^G_x\left(\frac{1}{\mu_x^G}+\max_{1\leq i\leq L+1}M_i\geq \lambda \sqrt {r(G)R_G(x,y)}\right),
\end{eqnarray*}
where $(M_i)_{i\geq 0}$ is the martingale defined at (\ref{mdef}). For $L\geq0$ and any $\theta\geq0$, we have that
\begin{eqnarray*}
\mathbf{P}^G_x\left(\max_{1\leq i\leq L+1}M_i\geq \lambda\right)&\leq& \mathbf{E}^G_x\left(e^{\theta M_{L+1}}\right)e^{-\theta\lambda}\\
&=& \mathbf{E}^G_x\left(e^{-\theta(\eta_i-1/\mu_x^G)}\right)^{L+1}e^{-\theta\lambda}\\
&\leq &e^{\frac{\theta^2}2\mathbf{E}^G_x\left(\eta_i^2\right)(L+1)}e^{-\theta\lambda}\\
&\leq &e^{{\theta^2}R_G(x,y)(L+1)/\mu_x^G-\theta\lambda},
\end{eqnarray*}
where we have applied Doob's sub-martingale inequality, $e^{-x}\leq 1-x+x^2$ for $x\geq 0$, and the second moment estimate of (\ref{2ndmom}), similarly to the proof of Theorem \ref{mainest}(a). We note that, because we are only seeking a one-sided bound, integrability is not an issue here, and hence no restrictions on $\theta\geq0$ are required. Optimising over $\theta\geq 0$ yields
\[\mathbf{P}^G_x\left(\max_{1\leq i\leq L+1}M_i\geq \lambda\right)\leq e^{-\lambda^2\mu_x^G/2R_G(x,y)(L+1)},\]
from which conclude that
\begin{equation}\label{yo}
\mathbf{P}^G_z\left(\max_{0\leq t\leq \tau_x(L+1)}\left(L^G_t(x)-L^G_t(y)\right)\geq \lambda \sqrt {r(G)R_G(x,y)}\right)\leq e^{-\lambda^2\mu_x^Gr(G)/8(L+1)}
\end{equation}
whenever $\lambda>2$.

Next, define the event
\[A(x,y,\lambda,L):=\left\{\max_{t\geq 0}\left(L\wedge \left(\frac{L^G_t(x)}{r(G)}\right)-L\wedge \left(\frac{L^G_t(y)}{r(G)}\right)\right)\geq \lambda\sqrt {\tilde{R}_G(x,y)}\right\}.\]
On this event, for some $t\geq 0$ we have that
\[\left(L\wedge \left(\frac{L^G_t(x)}{r(G)}\right)-L\wedge \left(\frac{L^G_t(y)}{r(G)}\right)\right)\geq \lambda\sqrt {\tilde{R}_G(x,y)},\]
and so, assuming $\lambda>0$ and $x\neq y$, it must also hold that $r(G)^{-1}L^G_t(y)<L$. In particular, this implies
\[L\wedge \left(\frac{L^G_t(x)}{r(G)}\right)\geq\lambda\sqrt {\tilde{R}_G(x,y)}+ \frac{L^G_t(y)}{r(G)}.\]
Setting $s=t\wedge \tau_x(\lfloor L\mu_x^Gr(G)\rfloor)$, we thus deduce that
\[\frac{L^G_s(x)}{r(G)}=L\wedge \left(\frac{L^G_t(x)}{r(G)}\right)\geq \lambda\sqrt {\tilde{R}_G(x,y)}+ \frac{L^G_t(y)}{r(G)}\geq \lambda\sqrt {\tilde{R}_G(x,y)}+ \frac{L^G_s(y)}{r(G)}.\]
Again applying the strong Markov property at $\tau_x(0)$, we have therefore shown that
\[\mathbf{P}_z^G\left(A(x,y,\lambda,L)\right)\leq \mathbf{P}_x^G\left(\max_{0\leq t\leq \tau_x(\lfloor L\mu_x^Gr(G)\rfloor)}\left(L^G_t(x)-L^G_t(y)\right)\geq \lambda \sqrt {r(G)R_G(x,y)}\right).\]
Recalling the bound at (\ref{yo}), this implies
\[\mathbf{P}_z^G\left(A(x,y,\lambda,L)\right)\leq  e^{-\lambda^2\mu_x^Gr(G)/8\lfloor L\mu_x^Gr(G)\rfloor}\leq e^{-\lambda^2/8 L}\]
for every $\lambda>2$ and $L\geq1$. The result follows.
\end{proof}

\section{A discrete version of Garsia's lemma}\label{garsec}

For our discrete version of Garsia's lemma, we continue to work in a general framework. In this section, though, we do not need to restrict our attention to the resistance metric, and instead consider an arbitrary metric $d_G$ on $V(G)$. We write $d_0(G):=\min_{x,y\in V(G):x\neq y}d_G(x,y)$ for the shortest distance between two distinct points,
\begin{equation}\label{ddiam}
d(G):=\max_{x,y\in V(G)}d_G(x,y)
\end{equation}
for the diameter of $V(G)$, and
\begin{equation}\label{dball}
B_d(x,r):=\left\{y:d_G(x,y)<r\right\}
\end{equation}
for the open balls with respect to this metric. To state our main result, we further suppose: $v:\mathbb{R}_+\rightarrow \mathbb{R}_+$ is a non-decreasing function; $p:\mathbb{R}_+\rightarrow \mathbb{R}_+$ is a non-decreasing function with $p(0)=0$; and $\psi:\mathbb{R}\rightarrow\mathbb{R}_+$ is symmetric, convex and satisfies $\psi(0)=1$ and $\lim_{x\rightarrow\infty}\psi(x)=\infty$.

{\propn\label{garsia} Suppose the measure $\mu^G$ satisfies
\begin{equation}\label{vollower}
\min_{x\in V(G)}\mu^G\left(B_d(x,r)\right)\geq v(r)
\end{equation}
for every $r=[d_0(G),d(G)]$. Given a function $f:V(G)\rightarrow\mathbb{R}$, define
\[\Gamma(f):= \sum_{x,y\in V(G)}\psi\left(\frac{f(x)-f(y)}{p(d_G(x,y))} \right)\mu_x^G\mu_y^G.\]
(By convention, we set $(f(x)-f(y))/p(d_G(x,y))=0$ when $x=y$.) It is then the case that
\begin{equation}\label{gsum}
\left|f(x)-f(y)\right|\leq2\sum_{i=1}^{\lfloor\log_2(d_G(x,y)/d_0(G))\rfloor+1}p(d_0(G)2^{i+1})\psi^{-1}\left(\frac{\Gamma(f)}{v(d_0(G)2^{i-1})^2}\right)
\end{equation}
for every $x,y\in V(G)$, where $\psi^{-1}(x):=\inf\{y\geq0:\psi(y)>x\}$.}
\begin{proof} Fix $x_0,y_0\in V(G)$, $x_0\neq y_0$. For $i=0,1,2,\dots$, define $A_i$ to be a set of the form $B_d(x_i,d_0(G)2^i)$ that contains $x_0$. Similarly, define $B_i$ to be a set of the form $B_d(y_i,d_0(G)2^i)$ that contains $y_0$. Note that $A_0=\{x_0\}$ and $B_0=\{y_0\}$. Moreover, we can choose the sets so that $A_n=B_n$, where $n:=\min\{i:d_0(G)2^i>d_G(x_0,y_0)\}=\lfloor\log_2(d_G(x_0,y_0)/d_0(G))\rfloor+1$.

For a set $A\subseteq V(G)$, if we define $f_A:=\frac{1}{\mu^G(A)}\sum_{x\in A} f(x)\mu_x^G$ to be the mean of $f$ on the set $A$, then, by applying the convexity of $\psi$, we deduce that, for $i=1,\dots,n$,
\begin{eqnarray*}
\psi\left(\frac{f_{A_i}-f_{A_{i-1}}}{p(d_0(G)2^{i+1})}\right)&\leq&\frac{1}{\mu^G(A_i)\mu^G(A_{i-1})}\sum_{x\in A_{i}}\sum_{y\in A_{i-1}}\psi\left(\frac{f(x)-f(y)}{p(d_0(G)2^{i+1})}\right)\mu^G_x\mu^G_y\\
&\leq &\frac{1}{\mu^G(A_{i-1})^2}\sum_{x\in A_{i}}\sum_{y\in A_{i-1}}\psi\left(\frac{f(x)-f(y)}{p(d_G(x,y))}\right)\mu^G_x\mu^G_y\\
&\leq & \frac{\Gamma(f)}{v(d_0(G)2^{i-1})^2}.
\end{eqnarray*}
In particular, this implies
\[\left|f_{A_i}-f_{A_{i-1}}\right|\leq p(d_0(G)2^{i+1})\psi^{-1}\left( \frac{\Gamma(f)}{v(d_0(G)2^{i-1})^2}\right).\]
Since $f_{A_0}=f(x_0)$, summing over $i$ gives
\[\left|f_{A_n}-f(x_0)\right|\leq \sum_{i=1}^n p(d_0(G)2^{i+1})\psi^{-1}\left( \frac{\Gamma(f)}{v(d_0(G)2^{i-1})^2}\right).\]
Repeating the argument for $y_0$ yields the desired result.
\end{proof}

{\rem\label{garsiarem} An elementary argument gives that the sum at (\ref{gsum}) can be bounded above by the integral
\[4\int_{d_0(G)}^{2d_G(x,y)}\frac{p(4s)}s\psi^{-1}\left(\frac{\Gamma(f)}{v(s/2)^2}\right)ds.\]
By extending the lower limit of integration to 0, one obtains an upper bound that does not depend on $d_0(G)$.}

\section{Local time continuity under UVD}\label{uvdsec}

In this section, we prove Theorem \ref{mainresult}. We start by combining Theorem \ref{mainest}(a) with the discrete version of Garsia's lemma derived in the previous section to establish a slightly weaker result. Although this has a worse power of the log term in the modulus of continuity than we will eventually obtain, it will allow us to uniformly control the maximum value of local time over the relevant time scales, as we do in the subsequent lemma.

{\lem\label{first} If $(G_i)_{i\in I}$ is a collection of graphs that satisfies UVD, then, for each $T>0$,
\[\lim_{\lambda\rightarrow\infty}\sup_{i\in I}\max_{z\in V(G_i)}\mathbf{P}^{G_i}_z\left(\max_{x,y\in V(G_i)}\max_{0\leq t\leq Tm(G_i)r(G_i)}\frac{r(G_i)^{-1}\left|L^{G_i}_{t}(x)-L_t^{G_i}(y)\right|}{{\tilde{R}_{G_i}(x,y)^{1/2}\left(1+\ln \tilde{R}_{G_i}(x,y)^{-1}\right)}}\geq \lambda\right)=0.\]}
\begin{proof} Set $d_{G_i}:=\tilde{R}_{G_i}$, $v_{G_i}(x):=v(r(G_i)x)$ (where $v$ is the function that appears in the definition of the UVD property), $p_{G_i}(x):=\sqrt{x}$ and $\psi_{G_i}(x):=e^{c|x|}$ for some $c$ which will later taken to be small. By the lower bound of UVD, we know that (\ref{vollower}) holds for each $i$ with this choice of $d_{G_i}$ and $v_{G_i}$. We therefore obtain from Proposition \ref{garsia} (and Remark \ref{garsiarem}) that
\[r(G_i)^{-1}\left|L^{G_i}_t(x)-L^{G_i}_t(y)\right|\leq
\frac{8}{c}\int_{0}^{2\tilde{R}_{G_i}(x,y)}\frac{1}{s^{1/2}}\ln_+\left(\frac{\Gamma\left({r(G_i)^{-1}L_t^{G_i}}\right)}{v(r(G_i)s/2)^2}\right)ds\]
for every $x,y\in V(G_i)$, $i\in I$, $t\geq　0$.

Assume now that $\Gamma({r(G_i)^{-1}L_t^{G_i}})\leq \lambda m(G_i)^2$ for some $\lambda\geq 1$. From the UVD property, we deduce
\[\Gamma\left({r(G_i)^{-1}L_t^{G_i}}\right)\leq \lambda m(G_i)^2\leq \lambda c_2^2 v(r(G_i))^2 \leq  \lambda c_2^2 c_3^4 s^{-2\log_2 c_3}v(r(G_i)s/2)^2,\]
where $c_3$ is the constant such that $v(2r)\leq c_3 v(r)$. Hence, setting $\tilde{\lambda}=\lambda^{1/2\log_2 c_3}$,
\begin{eqnarray*}
r(G_i)^{-1}\left|L^{G_i}_t(x)-L^{G_i}_t(y)\right|&\leq&
c_4\int_{0}^{2\tilde{R}_{G_i}(x,y)}\frac{1}{s^{1/2}}\ln_+\left(\frac{c_5\tilde{\lambda}}{s}\right)ds\\
\end{eqnarray*}
for every $x,y\in V(G_i)$, $i\in I$. This implies
\begin{eqnarray*}
r(G_i)^{-1}\left|L^{G_i}_t(x)-L^{G_i}_t(y)\right|&\leq& c_6 \sqrt{\tilde{R}_{G_i}(x,y)}\left(\ln(c_5\tilde{\lambda})+\ln \tilde{R}_{G_i}(x,y)^{-1}\right)\\
&\leq& c_7\left(1+\ln\lambda\right) \sqrt{\tilde{R}_{G_i}(x,y)}\left(1+\ln \tilde{R}_{G_i}(x,y)^{-1}\right)
\end{eqnarray*}
for every $x,y\in V(G_i)$, $i\in I$.

It follows from the conclusion of the previous paragraph that
\begin{eqnarray*}
\lefteqn{\sup_{i\in I}\max_{z\in V(G_i)}\mathbf{P}^{G_i}_z\left(\max_{x,y\in V(G_i)}\max_{0\leq t\leq Tm(G_i)r(G_i)}\frac{r(G_i)^{-1}\left|L^{G_i}_{t}(x)-L_t^{G_i}(y)\right|}{\tilde{R}_{G_i}(x,y)^{1/2}\left(1+\ln \tilde{R}_{G_i}(x,y)^{-1}\right)}\geq \lambda\right)}\\
&\leq &\sup_{i\in I}\max_{z\in V(G_i)}\mathbf{P}^{G_i}_z\left(\max_{0\leq t\leq Tm(G_i)r(G_i)}\Gamma\left({r(G_i)^{-1}L_t^{G_i}}\right)> \lambda' m(G_i)^2\right),\hspace{100pt}
\end{eqnarray*}
where $\lambda'$ is defined by $\lambda= c_7 (1+\ln \lambda')$. Hence, to complete the proof, it will be enough to show that
\begin{equation}\label{gambound}
\sup_{i\in I}\max_{z\in V(G_i)}m(G_i)^{-2}\mathbf{E}^{G_i}_z\left(\max_{0\leq t\leq Tm(G_i)r(G_i)}\Gamma\left({r(G_i)^{-1}L_t^{G_i}}\right)\right)<\infty.
\end{equation}

Now, by definition, we have that the left-hand side is bounded above by
\[\sup_{i\in I}\max_{x,y\in V(G_i)}\mathbf{E}^{G_i}_z\left(\exp\left\{\frac{c\max_{0\leq t\leq Tm(G_i)r(G_i)}r(G_i)^{-1}\left|L_t^{G_i}(x)-L_t^{G_i}(y)\right|}{\sqrt{\tilde{R}_{G_i}(x,y)}}\right\}\right).\]
Consequently, assuming that $c$ is chosen to be suitably small, the bound at (\ref{gambound}) can be deduced by applying Theorem \ref{mainest}(a).
\end{proof}

{\lem \label{suplem}If $(G_i)_{i\in I}$ is a collection of graphs that satisfies UVD, then, for each $T>0$,
\[\lim_{\lambda\rightarrow\infty}\sup_{i\in I}\max_{z\in V(G_i)}\mathbf{P}^{G_i}_z\left(\max_{x\in V(G_i)}r(G_i)^{-1}L^{G_i}_{Tm(G_i)r(G_i)}(x)\geq \lambda\right)=0.\]}
\begin{proof} First note that, for $\lambda\geq 1$,
$\mathbf{P}^{G}_x(r(G)^{-1}L^{G}_{Tm(G)r(G)}(x)> \lambda)$ is bounded above by $\mathbf{P}_x^G(\tau_{x}(\lfloor \lambda\mu_x^G r(G)\rfloor)< { Tm(G)r(G)})$.
Recalling the bound for $T_1$ that appeared in the proof of Theorem \ref{mainest}(a), it follows that
\[\mathbf{P}^{G_i}_x\left(r(G_i)^{-1}L^{G_i}_{Tm(G_i)r(G_i)}(x)> \lambda\right)\leq e^{-\mu_x^{G_i}\left(\frac{\lfloor \lambda\mu_x^{G_i} r(G_i)\rfloor}{\mu_x^{G_i}}-Tr(G_i)\right)^2/2c\lfloor \lambda\mu_x^{G_i} r(G_i)\rfloor r(G_i)},\]
and the upper bound here converges to 0 as $\lambda\rightarrow\infty$, uniformly in $x\in V(G_i)$, $i\in I$. The result follows by applying this convergence and Lemma \ref{first}, together with the observation that $\mathbf{P}^{G_i}_z(\max_{x\in V(G_i)}r(G_i)^{-1}L^{G_i}_{Tm(G_i)r(G_i)}(x)\geq \lambda)$ is bounded above by
\begin{eqnarray*}
\lefteqn{\mathbf{P}^{G_i}_z\left(r(G_i)^{-1}L^{G_i}_{Tm(G_i)r(G_i)}(z)\geq \lambda/2\right)}\\
&&+
\mathbf{P}^{G_i}_z\left(\max_{x,y\in V(G_i)}{r(G_i)^{-1}\left|L^{G_i}_{Tm(G_i)r(G_i)}(x)-L_{Tm(G_i)r(G_i)}^{G_i}(y)\right|}\geq \lambda/2\right).
\end{eqnarray*}
\end{proof}

{\lem\label{nearlem} If $(G_i)_{i\in I}$ is a collection of graphs that satisfies UVD, then, for each $L>0$,
\[\lim_{\lambda\rightarrow \infty}\sup_{i\in I}\max_{z\in V(G_i)}\mathbf{P}^{G_i}_z\left(\max_{x,y\in V(G_i)}\max_{t\geq0}\frac{\left|L\wedge \left(\frac{L^{G_i}_t(x)}{r(G_i)}\right)-L\wedge \left(\frac{L^{G_i}_t(y)}{r(G_i)}\right)\right|}{\sqrt{\tilde{R}_{G_i}(x,y)\left(1+\ln \tilde{R}_{G_i}(x,y)^{-1}\right)}}\geq \lambda\right)=0.\]}
\begin{proof} Since the proof of this is essentially the same as that of Lemma \ref{first} with the local times being replaced by the truncated local times, we omit the details. We merely note that to obtain the square root of the log term of the modulus of continuity, we take $\psi_{G_i}(x):=e^{cx^2}$ and estimate the expectation of $\max_{t\geq0}\Gamma ({L\wedge ({L^{G_i}_t}/{r(G_i)})})$
using the Gaussian bound of Theorem \ref{mainest}(b).
\end{proof}

\begin{proof}[Proof of Theorem \ref{mainresult}] Clearly, the probability we are trying to bound is no greater than
\begin{eqnarray*}
\lefteqn{\mathbf{P}^{G_i}_z\left(\max_{x,y\in V(G_i)}\max_{t\geq0}\frac{\left|L\wedge \left(\frac{L^{G_i}_t(x)}{r(G_i)}\right)-L\wedge \left(\frac{L^{G_i}_t(y)}{r(G_i)}\right)\right|}{\sqrt{\tilde{R}_{G_i}(x,y)\left(1+\ln \tilde{R}_{G_i}(x,y)^{-1}\right)}}\geq \lambda\right)}\\
&&+\mathbf{P}^{G_i}_z\left(\max_{x\in V(G_i)}r(G_i)^{-1}L^{G_i}_{Tm(G_i)r(G_i)}(x)\geq L\right).\hspace{60pt}
\end{eqnarray*}
Hence the result is an easy consequence of Lemmas \ref{suplem} and \ref{nearlem}.
\end{proof}

\section{Examples}\label{examsec}

In this section, we present some examples of collections of graphs for which UVD can be checked, and therefore to which Theorem \ref{mainresult} applies. Although in these examples we restrict our attention to collections of unweighted graphs (i.e.~those for which $\mu^{G}_{xy}=1$, for all $\{x.y\}\in E(G)$), we note that the assumption UVD is stable under perturbations that keep the weights uniformly bounded. In particular, the discussion would equally apply if we supposed $\mu^{G}_{xy}\in[c_1,c_2]$ for all $\{x.y\}\in E(G)$ (uniformly over the graphs in the collection), where $0<c_1\leq c_2<\infty$. We further note that for the majority of the graphs described explicitly in our examples, we have that for some $\alpha\geq 1$, $\beta\geq 2$,
\begin{equation}\label{volres}
\mu_G(B_d(x,r))\asymp r^\alpha,\qquad R_G(x,y)\asymp d_G(x,y)^{\beta-\alpha},
\end{equation}
for $x,y\in V(G)$, $r\in[r_0(G),r(G)]$ (again, uniformly over the collection), where $d_G$ is the usual shortest path graph distance, $B_d(x,r)$ is the corresponding ball (defined as at (\ref{dball})), and $\asymp$ means `bounded above and below by constant multiples of'. As a result, the conclusion of Theorem \ref{mainresult} can be written as
\begin{equation}\label{conc}
\lim_{\lambda\rightarrow\infty}\sup_{i\in I}\max_{z\in V(G_i)}\mathbf{P}^{G_i}_z\left(\max_{x,y\in V(G_i)}\max_{0\leq t\leq Td(G_i)^{\beta}}\frac{d(G_i)^{-(\beta-\alpha)}\left|L^{G_i}_{t}(x)-L_t^{G_i}(y)\right|}{\sqrt{\tilde{d}_{G_i}(x,y)^{\beta-\alpha}\left(1+\ln \tilde{d}_{G_i}(x,y)^{-1}\right)}}\geq \lambda\right)=0
\end{equation}
for each $T>0$, where $\tilde{d}_{G_i}(x,y):=d_{G_i}(x,y)/d(G_i)$ is the graph distance rescaled by the graph diameter (as defined at (\ref{ddiam})). Note that $\beta$ gives the relevant time-scaling exponent (cf.~the heat kernel estimates for infinite graphs of \cite{BCK}). We describe the extension of this local time continuity result to the infinite graph setting in the next section.

\subsection{One-dimensional graphs}

Consider $(G_i)_{i\in I}$ to be a collection of unweighted graphs for which there exists a finite constant $c$ such that $m(G_i)\leq c r(G_i)$ for all $i\in I$. Since the shortest path graph distance $d_{G_i}$ satisfies $d_{G_i}\geq R_{G_i}$, we immediately deduce that $\mu^{G_i}(B_{G_i}(x,r))\geq r$ for $r\in [r_0(G_i),r(G_i)]$, which confirms UVD holds in this case with $v(r)=r$. In particular, this class of examples covers collections of essentially one-dimensional graphs. For example, it includes the case when $G_i$ is the graph with vertices $\{0,\dots,i\}$ connected by nearest neighbour edges, $i\in \mathbb{N}$. (This latter example satisfies (\ref{volres}) uniformly over $(G_i)_{i\in \mathbb{N}}$ with $\alpha=1$, $\beta=2$.)

\subsection{Trees}

\begin{figure}[t]
\begin{center}
\scalebox{0.08}{\includegraphics{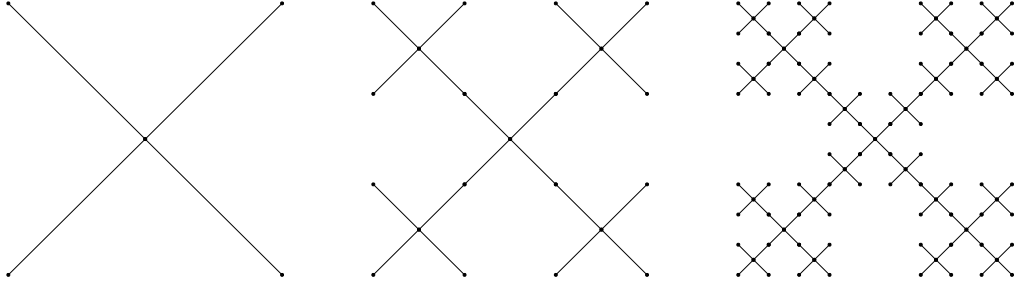}}
\end{center}
\vspace{-20pt}
\caption{The Vicsek set graphs $G_0$, $G_1$, $G_2$.}\label{vicsek}
\end{figure}

Suppose $(G_i)_{i\in I}$ is a collection of graph trees. Since in this case $d_{G_i}\equiv R_{G_i}$, it follows that we can replace the resistance metric by the shortest path metric in the UVD condition to be checked, and in the conclusion. In particular, if we have a family of trees with uniform polynomial volume growth of exponent $\alpha$ with respect to the graph distance (so that the left-hand estimate of (\ref{volres}) holds uniformly over $(G_i)_{i\in I}$), then UVD holds and Theorem \ref{mainresult} applies. (NB.~In this case, the right-hand estimate of (\ref{volres}) immediately holds with $\beta=\alpha+1$.) For instance, if we take the $G_i$ to be the $i$th level graph tree approximation of the Vicsek set -- the first three such graphs are shown in Figure \ref{vicsek}, then it is easy to check that we have the requisite polynomial volume growth with exponent $\alpha=\ln 5 / \ln 3$, and so we conclude (\ref{conc}) holds with this $\alpha$, $\beta=\alpha+1$ and $d(G_i)=2\times3^i$.

\subsection{Nested fractal graphs}\label{nestedsec}

The nested fractals were originally introduced in \cite{Lind}, and are a class of self-similar fractals that are finitely-ramified, embedded into Euclidean space and admit a high degree of symmetry. The volume and resistance growth of such fractals and associated graphs are well-understood, and so fit naturally into the framework of the present article. Although the discussion of this section would readily extend to any nested fractal, for simplicity of presentation we restrict ourselves to the graphs associated with the Sierpi\'nski gasket in two dimensions.

\begin{figure}[t]
\begin{center}
\scalebox{0.14}{\includegraphics{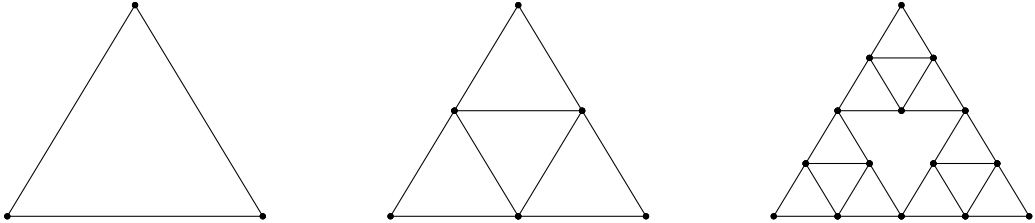}}
\end{center}
\vspace{-20pt}
\caption{The Sierpi\'nski gasket graphs $G_0$, $G_1$, $G_2$.}\label{sg}
\end{figure}

Let $V_0:=\{x_1,x_2,x_3\}\subseteq\mathbb{R}^2$ consist of the vertices of an equilateral triangle of side length 1. Write $\psi_i(x):=|x+x_i|/2$. Then there exists a unique compact set such $F$ that $F=\cup_{i=1}^3\psi_i(F)$; this is the Sierpi\'nski gasket. We define the associated Sierpi\'nski gasket graphs $(G_i)_{i\geq 0}$ by setting $V(G_i):=V_i$, where $V_i:=\cup_{i=1}^3\psi_i(V_{i-1})$
for $i\geq 1$, (note that $V_0$ was already defined,) and defining $E(G_i)$ to be the collection of pairs of elements of $V_i$ at a Euclidean distance $2^{-i}$ apart. (The first three graphs in this sequence are shown in Figure \ref{sg}.) For such graphs, it is easy to check that balls in the shortest path graph distance $d_{G_i}$ satisfy $c_1r^{d_f}\leq \mu^{G_i}(B_{d_{G_i}}(x,r))\leq c_2r^{d_f}$ for every $x\in V(G_i)$, $r\in [1,d(G_i)] $, $i\geq 0$, where $d_f:=\ln3/\ln 2$. Moreover, for the resistance metric $R_F$ on the limiting fractal, it is known that $c_3|x-y|^{d_w-d_f}\leq R_F(x,y)\leq c_4|x-y|^{d_w-d_f}$ for all $x,y\in F$, where $d_w=\ln 5/\ln 2$ (see \cite[(1.6.10)]{Str}, for example). From the standard construction of $R_F$ in terms of resistances on approximating subsets, it is possible to deduce that $R_{G_i}(x,y)=(5/3)^iR_F(x,y)$ for all $x,y\in V(G_i)$, $i\geq 0$. It is also straightforward to verify $c_52^i|x-y|\leq d_{G_i}(x,y)\leq c_62^i|x-y|$ for all $x,y\in V(G_i)$, $i\geq 0$. Hence it follows that
\[c_7d_{G_i}(x,y)^{d_w-d_f}\leq R_{G_i}(x,y)\leq c_8d_{G_i}(x,y)^{d_w-d_f}\]
for all $x,y\in V(G_i)$, $i\geq0$. (For nested fractals in general, a discussion of the connection between the various distances can be found in \cite[Remark 3.7]{FHK}.) Putting these estimates together, we deduce UVD holds for this example, and an application of Theorem \ref{mainresult} yields the following. We note that a similar modulus of continuity for the local times of the limiting diffusion was established in \cite[Theorem 1.11, see also the remark following its proof]{BP}.

{\thm \label{sgequi} If $(G_i)_{i\geq 0}$ is the sequence of Sierpi\'nski gasket graphs, then, for each $T>0$,
\[\lim_{\lambda\rightarrow\infty}\sup_{i\geq 0}\max_{z\in V(G_i)}\mathbf{P}^{G_i}_z\left(\max_{{x,y\in V(G_i)}}\max_{0\leq t\leq 5^iT}\frac{\left(\frac{3}{5}\right)^i\left|L^{G_i}_{t}(x)-L_t^{G_i}(y)\right|}{|x-y|^{\ln (5/3)/2\ln 2}(1+\ln |x-y|^{-1})^{1/2}}\geq \lambda\right)=0.\]}

\vspace{-30pt}
\subsection{Sierpi\'nski carpet graphs}\label{scsec}

\begin{figure}[t]
\begin{center}
\scalebox{0.08}{\includegraphics{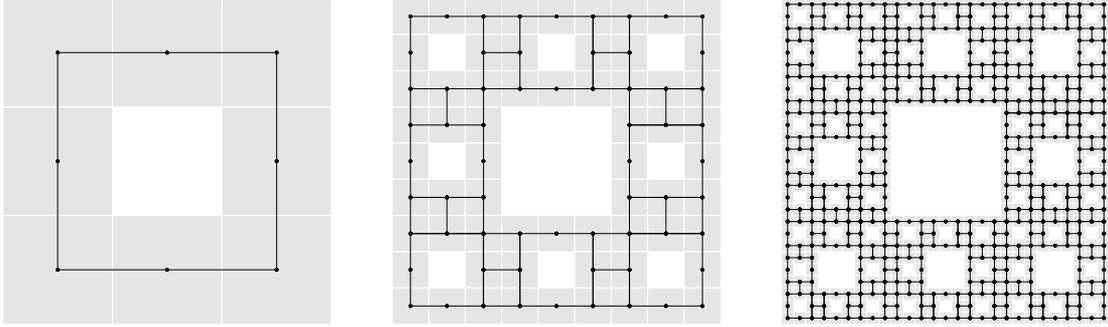}}
\end{center}
\vspace{-20pt}
\caption{The Sierpi\'nski carpet graphs $G_0$, $G_1$, $G_2$.}\label{sc}
\end{figure}

There are various definitions of generalised Sierpi\'nski carpets and associated graphs to which the following argument could be applied. Again, though, to avoid unnecessary complication, we take one representative example. Let $\{x_1,\dots,x_8\}\subseteq\mathbb{R}^2$ be the corners and edge-midpoints of the unit square $[0,1]^2$. Write $\psi_i(x):=|x+x_i|/3$. Then there exists a unique compact set such $F$ that $F=\cup_{i=1}^8\psi_i(F)$; this is the Sierpi\'nski carpet. We define the associated Sierpi\'nski carpets graphs $(G_i)_{i\geq 0}$ by first setting $V(G_i):=V_i$, where $V_0$ is the set consisting of the centres of the squares $(\psi_j([0,1]^2))_{j=1}^{8}$, and $V_i:=\cup_{i=1}^8\psi_i(V_{i-1})$ for $i\geq 1$. Moreover, we define $E(G_i)$ to be the collection of pairs of elements of $V_i$ at a Euclidean distance $3^{-(i+1)}$ apart. (The first three graphs in this sequence are shown in Figure \ref{sc}.) Let us also define an infinite version of the graphical Sierpi\'nski carpet $G$ by setting
$V(G):=\cup_{i=0}^\infty 3^{i+1}V_i$, and defining $E(G)$ to be the collection of pairs of elements of $V(G)$ a unit distance apart.

Now, for the infinite graphical Sierpi\'nski carpet, it is a consequence of results in \cite{BBRes, BCK} that
\begin{equation}\label{igsc}
c_1d_{G}(x,y)^{d_w-d_f}\leq R_{G}(x,y)\leq c_2d_{G}(x,y)^{d_w-d_f}
\end{equation}
for every $x,y\in V(G)$, where $R_G$ and $d_G$ are the resistance and shortest path metric, respectively, $d_f:=\ln 8/\ln 3$ and $d_w:= \ln (8\rho)/\ln 3$ for some $\rho>1$. Hence, if we consider $G_i$ with `wired' boundary conditions, by which we mean we identify all the vertices on the outer edge of the graph to obtain a new graph $G_i^{\rm w}$, then a straightforward application of Rayleigh's monotonicity law (see, for example, \cite[Theorem 9.12]{LPW}) allows us to deduce that $R_{G_i^{\rm w}}(x,y)\leq c_3d_{G_i^{\rm w}}(x,y)^{d_w-d_f}$ for all $x,y\in V(G_i^{\rm w})$, $i\geq 0$. It is also an elementary exercise to check that
$\mu^{G_i^{\rm w}}(B_{d_{G_i^{\rm w}}}(x,r))\geq \mu^{G_i}(B_{d_{G_i}}(x,r))\geq c_4r^{d_f}$ for every $x\in V(G_i^{\rm w})$, $r\in [1,d(G_i^{\rm w})] $, $i\geq 0$. (If $x$ is the boundary vertex in $G_i^{\rm w}$, then we can take any $x$ on the boundary in the middle expression above.) This confirms that the first part of the UVD condition holds with $v(r)=r^{d_f/(d_f-d_w)}$.

\begin{figure}[t]
\begin{center}
\scalebox{1}{\includegraphics{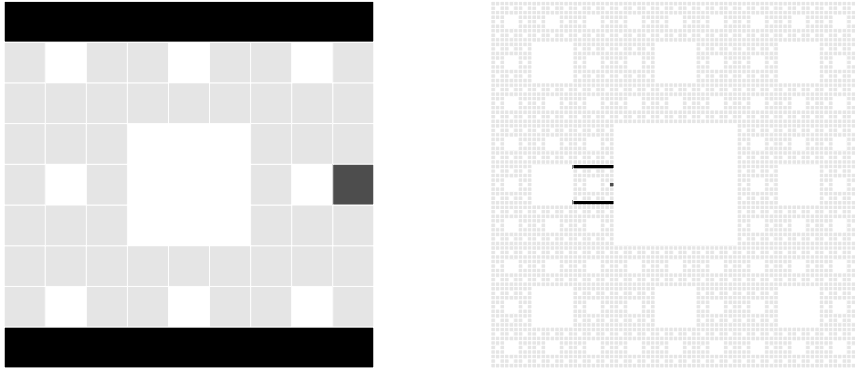}}
\end{center}
\vspace{-20pt}
\caption{The left-hand figure shows $A$ (black) and $B$ (dark grey). The right-hand figure shows $A'$ and $B'$.}\label{sc2}
\end{figure}

For the second part of the UVD condition, let us start by defining $A$ to be the union of two rectangles of height $1/9$ and width 1, one at the top and one at the bottom of the unit square $[0,1]^2$. Moreover, define $B$ to be the square of side $1/9$ located on the middle of the right-hand side of $[0,1]^2$. (See left-hand side of Figure \ref{sc2}.) From \cite[Theorem 6.1]{McG}, it follows that there exists a constant $c_5$ such that
\begin{equation}\label{reslower}
R_{G_i}(A\cap V_i,B\cap V_i)\geq c_5 \rho^i
\end{equation}
for every $i\geq 1$. (The graphs considered in \cite{McG} have larger vertex sets than ours, but it is easy to see that the resistance in the two settings is comparable.) Next, let $A'$ and $B'$ be the image of $A$ and $B$ under the map that takes the unit square $[0,1]^2$ to $[2/9,1/3]\times [4/9,5/9]$ (see the right-hand side of Figure \ref{sc2}). By again applying Rayleigh's monotonicity law, it follows from  (\ref{reslower}) that if $x\in B'\cap V(G_i^{\rm w})$ and $y\in V(G_i^{\rm w})\backslash [2/9,1/3]\times [4/9,5/9]$ (where $i\geq 3$), then $R_{G_i^{\rm w}}(x,y)\geq R_{G_i^{\rm w}}(A',B')\geq c_6 \rho^i$. In particular, this implies that $r(G_i^{\rm w})\geq c_7 \rho^i$, and we conclude that $m(G_i^{\rm w})\leq c_7 8^i\leq c_8r(G_i^{\rm w})^{d_f/(d_f-d_w)}$, as desired. Thus, since the graph distance $d_{G_i}$ is comparable to the wired Euclidean distance $|\cdot-\cdot|_{\rm w}$ (i.e.~the quotient of the usual Euclidean distance on $[0,1]^2$ under the identification of the boundary) multiplied by $3^i$, we obtain the following. We expect that by developing further the techniques of \cite{McG} (see also \cite{BBRes}), it should be possible to verify the corresponding result for the unwired graphs $(G_i)_{i\geq0}$.

{\thm If $(G_i^{\rm w})_{i\geq 0}$ is the sequence of wired Sierpi\'nski carpet graphs, then, for each $T>0$,
\[\lim_{\lambda\rightarrow\infty}\sup_{i\geq 0}\max_{z\in V(G_i^{\rm w})}\mathbf{P}^{G_i^{\rm w}}_z\left(\max_{{x,y\in V(G_i^{\rm w})}}\max_{0\leq t\leq (8\rho)^iT}\frac{\rho^{-i}\left|L^{G_i^{\rm w}}_{t}(x)-L_t^{G_i^{\rm w}}(y)\right|}{|x-y|_{\rm w}^{\ln (\rho)/2\ln 3}(1+\ln |x-y|_{\rm w}^{-1})^{1/2}}\geq \lambda\right)=0.\]}

\section{Infinite graphs}\label{infsec}

In this section, we consider the application of the techniques developed in this article to a class of infinite graphs. In particular, we suppose $G=(V(G),E(G))$ is an infinite, locally finite, connected graph, with weights $\mu^G$ and distinguished vertex $0\in V(G)$. We assume that this satisfies (\ref{volres}) for some $\alpha<\beta$ (where $\alpha\geq 1$, $\beta\geq 2$). We note that, as is discussed following \cite[Definition 1.2]{BCK} (see also \cite[Proposition 3.5]{BCK}), these conditions imply that the random walk $X^G$ is recurrent. As a consequence, the identity at (\ref{returnprob}) still holds in this setting (see \cite[Lemma 2.48]{Barlowrg}, for example). This will be useful in proving the following adaptation of Theorem \ref{mainresult} to the present setting, which is the main result of this section. Since the proof is in many aspects similar to that of Theorem \ref{mainresult}, we will be brief with the details. At the end of the section, we describe a particular application to the infinite graphical Sierpi\'{n}ski carpet.

{\thm \label{infthm} If $G$ is an infinite graph satisfying (\ref{volres}) for some $\alpha\geq 1$, $\beta\geq 2$ such that $\alpha<\beta$, then, for each $T>0$,
\[\lim_{\lambda\rightarrow\infty}\sup_{i\geq 1}\mathbf{P}^{G}_0\left(\max_{x,y\in V(G)}\max_{0\leq t\leq Ti^\beta}\frac{i^{-(\beta-\alpha)}\left|L^{G}_{t}(x)-L_t^{G}(y)\right|}{\sqrt{{R}_{G}^{(i)}(x,y)\left(1+\ln_+ \left({R}^{(i)}_{G}(x,y)\right)^{-1}\right)}}\geq \lambda\right)=0,\]
where ${R}^{(i)}_{G}(x,y):=i^{-(\beta-\alpha)}{R}_{G}(x,y)$.}
\bigskip

To prove the above result, we start with a lemma that controls the rate of growth of the local times at a given vertex. Throughout this section, we suppose $d_G$ is the usual shortest path metric on $G$, and denote balls with respect to this metric by $B_d(x,r)$.

{\lem \label{infsup} If $G$ is an infinite graph satisfying the conditions of Theorem \ref{infthm}, then, for each $T>0$, there exist constants $c_1$ and $c_2$ such that
\[\sup_{i\geq 1}\sup_{x\in V(G)}\mathbf{P}^{G}_x\left(i^{-(\beta-\alpha)}L^{G}_{Ti^\beta}(x)\geq \lambda\right)\leq c_1e^{-c_2\lambda}\]
for every $\lambda\geq 0$.}
\begin{proof} We will first show the existence of constants $\varepsilon_1,\varepsilon_2>0$ such that
\begin{equation}\label{returncond}
\mathbf{P}^{G}_x\left(\tau_x^+\geq \varepsilon_1i^\beta\right)\geq \frac{\varepsilon_2}{\mu_x^Gi^{\beta-\alpha}}
\end{equation}
for all $x\in V(G)$, $i\geq 1$. Given $x\in V(G)$, let $y\in V(G)$ be such that $d_G(x,y)=i$. We then have that
\begin{eqnarray*}
\mathbf{P}^{G}_x\left(\tau_x^+\geq \varepsilon_1i^\beta\right)&\geq &\mathbf{P}_x^{G}\left(\tau_y<\tau_x^+\right)\mathbf{P}_y^{G}\left(\tau_{B_d(y,i)^c}>\varepsilon_1 i^\beta\right)\\
&\geq& \frac{1}{\mu_x^GR_G(x,y)}\left(1-c_1e^{-c_2\varepsilon_1^{-\frac{1}{\beta-1}}}\right)\\
&\geq &  \frac{c_3}{\mu_x^Gi^{\beta-\alpha}},
\end{eqnarray*}
where $\tau_{B_d(y,i)^c}$ is the exit time of the ball $B_d(y,i)$. Note that we have applied (\ref{returnprob}) and \cite[Proposition 3.4 and Lemma 3.7]{BCK} to deduce the second inequality, and (\ref{volres}) to obtain the third. This confirms the desired bound.

Now, define $(\tau_x(i))_{i\geq 0}$ as in the proof of Theorem \ref{mainest}(a). It is then the case that
\begin{eqnarray*}
\mathbf{P}^{G}_x\left(i^{-(\beta-\alpha)}L^{G}_{Ti^\beta}(x)>\lambda\right)&\leq&
\mathbf{P}^{G}_x\left(\tau(\lfloor\lambda \mu_x^Gi^{\beta-\alpha} \rfloor)<Ti^\beta\right)\\
&\leq &\mathbf{P}_x^G\left(\sum_{j=0}^{\lfloor\lambda \mu_x^Gi^{\beta-\alpha} \rfloor-1}\mathbf{1}_{\{\tau_{x}(j+1)-\tau_{x}(j)\geq \varepsilon_1i^\beta\}}< T/\varepsilon_1 \right)\\
&=&\mathbf{P}\left({\rm Bin} \left(\lfloor\lambda \mu_x^Gi^{\beta-\alpha} \rfloor,\mathbf{P}^G_x(\tau_{x}^+\geq \varepsilon_1 i^\beta)\right)<T/\varepsilon_1 \right)\\
&\leq &\mathbf{P}\left({\rm Bin} \left(\lfloor\lambda \mu_x^Gi^{\beta-\alpha} \rfloor,\varepsilon_2/\mu_x^Gi^{\beta-\alpha}\right)< T/\varepsilon_1 \right),
\end{eqnarray*}
where we denote by ${\rm Bin}(n,p)$ a binomial random variable with parameters $n$ and $p$, built on a probability space with probability measure $\mathbf{P}$. Note that the final inequality is a consequence of (\ref{returncond}). Hence,
\begin{eqnarray*}
\mathbf{P}^{G}_x\left(i^{-(\beta-\alpha)}L^{G}_{Ti^\beta}(x)>\lambda\right)&\leq&e^{T/\varepsilon_1}
\mathbf{E}^{G}_x\left(e^{-{\rm Bin} \left(\lfloor\lambda \mu_x^Gi^{\beta-\alpha} \rfloor,\varepsilon_2/\mu_x^Gi^{\beta-\alpha}\right)}\right)\\
&\leq &e^{T/\varepsilon_1}e^{-(1-e^{-1})\lfloor\lambda \mu_x^Gi^{\beta-\alpha} \rfloor\varepsilon_2/\mu_x^Gi^{\beta-\alpha}}.
\end{eqnarray*}
Next, note that if $y\in V(G)$ is such that $d_G(x,y)=1$, then, by applying (\ref{returnprob}) and (\ref{volres}), we have that $1\leq \mathbf{P}_{x}(\tau_y<\tau_x^+)^{-1}=\mu_x^GR_G(x,y)\leq c_4\mu_x^Gd_G(x,y)^{\beta-\alpha}=c_4\mu_x^G$. In conjunction with the above inequality, it follows that, for $\lambda\geq 2c_4$, $\mathbf{P}^{G}_x(i^{-(\beta-\alpha)}L^{G}_{Ti^\beta}(x)>\lambda)\leq c_5e^{-c_6\lambda}$, and the result follows.
\end{proof}

The following result is a version of Theorem \ref{mainest} for infinite graphs. Since on replacing $r(G)$ by $i^{\beta-\alpha}$, and $m(G)$ by $i^{\alpha}$, the proof of the result is almost identical to that of Theorem \ref{mainest}, we omit it. (The one other change that is required is the use of Lemma \ref{infsup} to bound the term corresponding to $T_1$ in the proof of part (a).)

{\lem \label{infmainest} Suppose $G$ is an infinite graph satisfying the conditions of Theorem \ref{infthm}.\\
(a) For each $\kappa, T>0$, there exist constants $c_1$ and $c_2$ such that
\[\sup_{i\geq 1}\max_{\substack{x,y,z\in V(G):\\d_G(x,y)\leq \kappa i}}\mathbf{P}^G_z\left(\max_{0\leq t\leq Ti^\beta}i^{-(\beta-\alpha)}\left|L^G_t(x)-L^G_t(y)\right|\geq \lambda\sqrt {{R}^{(i)}_G(x,y)}\right)\leq c_1e^{-c_2 \lambda}\]
for every $\lambda\geq 0$.\\
(b) It holds that
\[\sup_{i\geq 1}\max_{x,y,z\in V(G)}\mathbf{P}^G_z\left(\max_{t\geq 0}\left|L\wedge \left(\frac{L^G_t(x)}{i^{\beta-\alpha}}\right)-L\wedge \left(\frac{L^G_t(y)}{i^{\beta-\alpha}}\right)\right|\geq \lambda\sqrt {{R}^{(i)}_G(x,y)}\right)\leq 2e^{\frac12-\frac{\lambda^2}{8L}}\]
for every $\lambda\geq 0$ and $L\geq 1$.}
\bigskip

We are now in a position to prove the main result of the section.

\begin{proof}[Proof of Theorem \ref{infthm}] Given $\kappa\geq 6$, let $G_i$ be the graph with vertex set $V(G_i):=B_d(0,\kappa i)$ and edge set $E(G_i):=\{x,y\in E(G):\: x,y\in B_d(0,\kappa i)\}$. From (\ref{volres}), it is possible to check that there exists a constant $c_1$ such that
\begin{equation}\label{mmm}
\mu^{G}\left(B_d(x,r)\cap B_d(0,\kappa i)\right)\geq c_1r^{\alpha}
\end{equation}
for every $x\in V(G_i)$, $r\in[1,2\kappa i]$, $i\geq 1$. Indeed, for balls such that $B_d(x,r)\subseteq B_d(0,\kappa i)$ (which includes the case $r=1$) this is obvious. Otherwise, $x\in B_d(0,(\kappa i-r)_+)^c$. If we further suppose $2\leq r\leq 2\kappa i/3 $, then $1\leq \lfloor r/2\rfloor\leq \kappa i-r$, and it is possible to select $y\in V(G_i)$ to be a point on a shortest path from $0$ to $x$ such that $d_G(x,y)= \lfloor r/2\rfloor$. For this $y$, we have $\mu^G(B_d(x,r)\cap B_d(0,\kappa i))\geq\mu^G(B_d(y, \lfloor r/2\rfloor))\geq c_2 r^\alpha$. On the other hand, assume  $2\kappa i/3 \leq r\leq 2\kappa i$. Let  $y\in V(G_i)$ to be a point on a shortest path from $0$ to $x$ such that $d_G(x,y)= \lfloor d_G(0,x)/2\rfloor$. Then $\mu^G(B_d(x,r)\cap B_d(0,\kappa i))\geq\mu^G(B_d(y,\kappa i/6))\geq c_3 (\kappa i)^\alpha\geq c_4r^\alpha$. This confirms (\ref{mmm}), and thus we deduce that the measures $\mu^{G_i}:=\mu^G(\cdot\cap B_d(0,\kappa i))$ satisfy (\ref{vollower}) uniformly in $i\geq 1$ for $v(r)=c_1 r^{\alpha}$.

Next, from Lemma \ref{infmainest}(a), we deduce that, for each $\kappa, T>0$, there exists a constant $c_5$ such that
\[\sup_{i\geq 1}i^{-2\alpha}\mathbf{E}^G_0\left(\max_{0\leq t\leq Ti^\beta}\sum_{x,y\in B_d(0,\kappa i)}e^{c_5i^{-(\beta-\alpha)}|L_t(x)-L_t(y)|/\sqrt {{R}^{(i)}_G(x,y)}}\mu_x^G\mu_y^G\right)<\infty.\]
Hence, setting $d_{G_i}:={R}_{G}^{(i)}$, $v_{G_i}(x):=i^\alpha x^{\alpha/(\beta-\alpha)}$, $p_{G_i}(x):=\sqrt{x}$ and $\psi_{G_i}(x):=e^{c|x|}$ for suitably small $c$, and applying the volume bound of the previous paragraph, one can proceed as in the proof of Lemma \ref{first} to show that, for any $\kappa, T>0$,
\[\lim_{\lambda\rightarrow\infty}\sup_{i\geq 1}\mathbf{P}^{G}_0\left(\max_{x,y\in B_d(0,\kappa i)}\max_{0\leq t\leq Ti^\beta}\frac{i^{-(\beta-\alpha)}\left|L^{G}_{t}(x)-L_t^{G}(y)\right|}{{{R}^{(i)}_{G}(x,y)^{1/2}\left(1+\ln_+ {R}^{(i)}_{G}(x,y)^{-1}\right)}}\geq \lambda\right)=0.\]
Together with the conclusion of Lemma \ref{infsup}, this implies that, for each $\kappa, T>0$,
\begin{equation}\label{gbound1}
\lim_{\lambda\rightarrow\infty}\sup_{i\geq 1}\mathbf{P}^{G}_0\left(\max_{x\in B_d(0,\kappa i)}i^{-(\beta-\alpha)}L^{G}_{Ti^\beta}(x)\geq \lambda\right)=0
\end{equation}
(cf.~the proof of Lemma \ref{suplem}). Moreover, Lemma \ref{infmainest}(b) implies that for each $\kappa>0$, $L\geq1$, there exists a constant $c_6$ such that
\[\sup_{i\geq 1}i^{-2\alpha}\mathbf{E}^G_0\left(\max_{t\geq 0}\sum_{x,y\in B_d(0,\kappa i)}e^{c_6\left|L\wedge \left(\frac{L^G_t(x)}{i^{\beta-\alpha}}\right)-L\wedge \left(\frac{L^G_t(y)}{i^{\beta-\alpha}}\right)\right|^2/{{R}^{(i)}_G(x,y)}}\mu_x^G\mu_y^G\right)<\infty.\]
Thus, as in Lemma \ref{nearlem}, taking $\psi_{G_i}(x):=e^{cx^2}$ in a similar application of Garsia's lemma yields
\begin{equation}\label{gbound2}
\lim_{\lambda\rightarrow \infty}\sup_{i\geq 1}\mathbf{P}^{G}_0\left(\max_{x,y\in B_d(0,\kappa i)}\max_{t\geq0}\frac{\left|L\wedge \left(\frac{L^{G}_t(x)}{i^{\beta-\alpha}}\right)-L\wedge \left(\frac{L^{G}_t(y)}{i^{\beta-\alpha}}\right)\right|}{\sqrt{{R}^{(i)}_{G}(x,y)\left(1+\ln_+ {R}^{(i)}_{G}(x,y)^{-1}\right)}}\geq \lambda\right)=0.
\end{equation}

Putting (\ref{gbound1}) and (\ref{gbound2}) together, we can emulate the proof of Theorem \ref{mainresult} to obtain
\begin{equation}\label{gbound3}
\lim_{\lambda\rightarrow\infty}\sup_{i\geq 1}\mathbf{P}^{G}_0\left(\max_{x,y\in B_d(0,\kappa i)}\max_{0\leq t\leq Ti^\beta}\frac{i^{-(\beta-\alpha)}\left|L^{G}_{t}(x)-L_t^{G}(y)\right|}{
\sqrt{{R}^{(i)}_{G}(x,y)\left(1+\ln_+ {R}^{(i)}_{G}(x,y)^{-1}\right)}}\geq \lambda\right)=0.
\end{equation}
Finally, we know from \cite[Proposition 3.4 and Lemma 3.7]{BCK} that if $\tau_{B_d(0,\kappa i)^c}$ is the exit time of $X^G$ from $B_d(0,\kappa i)$, then $\lim_{\kappa\rightarrow\infty}\sup_{i\geq 1}\mathbf{P}_0^G(\tau_{B_d(0,\kappa i)^c}\leq T i^\beta)=0$. The result readily follows by applying this together with (\ref{gbound3}).
\end{proof}

A simple application of Theorem \ref{infthm} yields the following corollary for the infinite graphical Sierpi\'{n}ski carpet introduced in Section \ref{scsec}. The key estimate verifying (\ref{volres}) is stated at (\ref{igsc}). (The relevant volume bound is easy to check.) In stating the result, we write the infinite carpet $F_{\infty}:=\cup_{i=0}^\infty 3^i F$, where $F$ is the Sierpi\'{n}ski carpet defined in Section \ref{scsec}. We extend the definition of discrete local times to this set in such a way that for each $x\in F_\infty$, $L^G_{t}(x)$ is bounded below (above) by the minimum (maximum) of $L^G_t$ over the graph vertices contained in the same and adjacent unit squares, and also $(L^G_{t}(x))_{x\in {F_\infty}}$ is continuous. We then extend to all times by linear interpolation. As is explained in the next section, were it known that the random walks converged under rescaling to a diffusion with jointly continuous local times, then this result would be enough to confirm that the local times of the random walk also converged under suitable rescaling.

{\cor For the infinite graphical Sierpi\'{n}ski carpet, for each $t\geq 0$, the laws of
\[\left(i^{-(d_w-d_f)}L^G_{i^{d_w} t}(ix)\right)_{x\in F_{\infty}},\]
$i=1,2,\dots$, form a tight sequence of probability measures on $C(F_{\infty},\mathbb{R})$, where $d_f:=\ln 8/\ln 3$ and $d_w:= \ln (8\rho)/\ln 3$.}

\section{Local time and cover time scaling}\label{scalingsec}

In this section, we consider the implications of local time equicontinuity for sequences of graphs for which the associated random walks admit a scaling limit. As in Section \ref{nestedsec}, for brevity we restrict ourselves to the unweighted Sierpi\'nski gasket graphs. It should be noted, however, that the arguments below are relatively generic, and should be transferable to many other models once the relevant inputs are established. Indeed, this is the reason why, despite it being possible to prove a stronger result for cover times than the one we derive below using a simple time-change argument in the particular case of nested fractal graphs, we believe the techniques developed here are still of interest (see Remarks \ref{f1} and \ref{f2} for further discussion on this point).

Let $(G_i)_{i\geq 0}$ be the sequence of Sierpi\'nski gasket graphs of Section \ref{nestedsec}, and $F$ be the limiting Sierpi\'nski gasket into which these are embedded. By \cite{BP,Goldstein,Kus}, we know that if the associated random walks $X^{G_i}$ are started from $x_i\in V(G_i)$, where $x_i\rightarrow x\in F$, then
\begin{equation}\label{scaling}
\left(X^{G_i}_{5^it}\right)_{t\geq 0}\rightarrow\left(X^{F}_{t}\right)_{t\geq 0}
\end{equation}
in distribution in $C(\mathbb{R}_+,F)$, where $X^F$ is Brownian motion on the Sierpi\'nski gasket started from $x$. (We suppose discrete time processes are extended to elements of $C(\mathbb{R}_+,F)$ by linear interpolation.)

In \cite{BP}, it was shown that the Brownian motion $X^F$ admits local times $(L_t^F(x))_{x\in F, t\geq 0}$ that, almost-surely, are jointly continuous in $(x,t)$ and satisfy the occupation density formula:
\[\int_Ff(x)L_t^F(x)\mu^F(dx) =\int_0^tf(X_s^F)ds,\]
for any continuous function $f:F\rightarrow \mathbb{R}$ and $t\geq 0$, where $\mu^F$ is the $(\ln3/\ln2)$-dimensional Hausdorff measure on $F$, normalised to be a probability measure. For $t\in\mathbb{N}$, we similarly have
\[\int_{F}f(x)L_t^{G_i}(x){\mu^{G_i}(dx)}=\sum_{j=0}^{t-1}f(X^{G_i}_j).\]
Hence, by applying the random walk scaling limit of (\ref{scaling}), it is possible to check that, for each continuous $f:F\rightarrow \mathbb{R}$ and $t\geq 0$,
\[5^{-i}\int_{F}f(x)L_{5^it}^{G_i}(x){\mu^{G_i}(dx)}\rightarrow \int_Ff(x)L_t^F(x)\mu^F(dx)\]
in distribution. (Note that, for this statement to make sense, we suppose that the definition of the discrete local time processes it extended to all times by linear interpolation at each vertex.) By \cite[Theorem 16.16]{Kall}, this is enough to imply that, for each $t\geq 0$,
\begin{equation}\label{msense}
5^{-i}L_{5^it}^{G_i}(x){\mu^{G_i}(dx)}\rightarrow L_t^{F}(x){\mu^{F}(dx)}
\end{equation}
in distribution in the topology of weak convergence of Borel measures on $F$. (We view $\mu^{G_i}$ as an atomic measure on $F$ in the obvious way.)

Now, for each $t$, we extend $(L^{G_i}_t(x))_{x\in V(G_i)}$ to a continuous function on $F$ by setting
\[L_t^{G_i}(x)=\frac{\sum_{k=1}^3 |x-x_k|^{-1} L_t^{G_i}(x_k)}{\sum_{k=1}^3 |x-x_k|^{-1}}\]
when $x$ is contained in the $i$th level triangle with vertices $x_1$, $x_2$, $x_3$. Given the equicontinuity result of Theorem \ref{sgequi} and uniform boundedness of Lemma \ref{suplem}, we can apply the Arzela-Ascoli theorem to deduce that the laws of $((3/5)^iL^{G_i}_{5^it}(x))_{x\in F}$ form a tight sequence of probability measures on $C(F,\mathbb{R}_+)$. In particular, the sequence $((3/5)^iL^{G_i}_{5^it}(x))_{x\in F}$ admits a distributionally convergent subsequence. Suppose that we have such a subsequence $((3/5)^{i_j}L^{G_{i_j}}_{5^{i_j}t}(x))_{x\in F}$, and $(\ell(x))_{x\in F}$ is the distributional limit in $C(F,\mathbb{R})$. Since $\mu^G/(6\times3^i)\rightarrow \mu^F$, it is an easy application of the continuous mapping theorem to deduce from this that
$5^{-{i_j}}L_{5^{i_j}t}^{G_{i_j}}(x){\mu^{G_{i_j}}(dx)}\rightarrow 6\ell(x){\mu^{F}(dx)}$ in distribution in the topology of weak convergence of Borel measures on $F$. In conjunction with (\ref{msense}) and the almost-sure continuity of $(L^F_t(x))_{x\in F}$, it follows that $(6\ell(x))_{x\in F}$ is equal to $(L^F_t(x))_{x\in F}$ in distribution. Since this conclusion is independent of the particular subsequence chosen, we obtain that, for each $t\geq 0$,
\[\left(6\left(\frac{3}{5}\right)^{i}L^{G_{i}}_{5^it}(x)\right)_{x\in F}\rightarrow (L^F_t(x))_{x\in F}\]
in distribution in $C(F,\mathbb{R}_+)$. Given that the convergence of the rescaled $X^{G_i}$ to $X^F$ holds in the uniform topology over compact time intervals, this result is readily extended to hold simultaneously over a finite collection of times $0\leq t_1\leq\dots \leq t_k$. Moreover, because local times are increasing in $t$ and the limit is continuous in the temporal variable, an elementary argument allows us to deduce the convergence is also uniform over time (cf.~the proof of Dini's theorem). In particular, by following these steps, we obtain the following result.

{\thm \label{localscaling} Let $(G_i)_{i\geq 0}$ be the sequence of Sierpi\'nski gasket graphs of Section \ref{nestedsec}. If the associated random walks $X^{G_i}$ are started from $x_i\in V(G_i)$, where $x_i\rightarrow x$, then
\[\left(6\left(\frac{3}{5}\right)^iL^{G_i}_{5^i t}(x)\right)_{x\in F, t\geq 0}\rightarrow
\left(L^{F}_{t}(x)\right)_{x\in F, t\geq 0}\]
in distribution in $C(F\times \mathbb{R}_+,\mathbb{R})$, where $(L_t^F(x))_{x\in F, t\geq 0}$ are the local times for the Brownian motion $X^F$ on the Sierpi\'nski gasket $F$ started from $x$.}

{\rem \label{f1} We now discuss a simpler proof of the corresponding result for the continuous time version of the random walk, similar to the proof of \cite[Theorem 7.22]{Barlow}. If we define $A^i_t:=\int_FL^F_t(x)\mu^{G_i}(dx)/m(G_i)$ and $\tau_i(t):=\inf\{s:A_s^i> t\}$, then $(X^F_{\tau_i(t)})_{t\geq 1}$ gives the continuous time random walk on $G_i$ with exponential mean $5^{-i}$ holding times. Moreover, similarly to the argument of \cite[Lemma 3.4]{Croy1}, one can check that the local times of this process with respect to $\mu^{G_i}/m(G_i)$ are given by $(L_{\tau_i(t)}^F(x))_{x\in V(G_i),t\geq0}$. Now, since $\mu^{G_i}/m(G_i)\rightarrow \mu^F$, the continuity of the local times $L^F$ imply that $(\tau_{i}{(t)})_{t\geq0}\rightarrow (t)_{t\geq0}$ almost-surely. Thus a simple reparametrisation yields that if $(\tilde{L}_{t}^{G_i}(x))_{x\in V(G_i),t\geq0}$ are the local times of the continuous time simple random walk on $G_i$ with exponential mean one holding times, with respect to the measure $\mu^{G_i}$, then $(6(3/5)^i\tilde{L}_{5^it}^{G_i}(x))_{x\in F,t\geq0}\rightarrow
({L}_{t}^{F}(x))_{x\in F,t\geq0}$ in distribution in $C(F\times \mathbb{R}_+,\mathbb{R})$, where the local times on the discrete spaces are suitably extended to take values in this space. We note, however, that this result would not transfer to the discrete time case without the use of some form of equicontinuity property, such as the one we have proved in this article. Moreover, the construction of the time-changed processes that the proof depends on is quite specific to the particular situation, and would not readily transfer to other settings, such as the Sierpi\'nski carpet.}
\bigskip

To conclude the article, we show that, as a consequence of this local time convergence, we are able to deduce the asymptotic behaviour of the cover times of the Sierpi\'nski gasket graphs in the sequence. To this end, for a random walk $X^G$ on a graph $G$, we define
\begin{equation}\label{tdef}
\tau_{\rm cov}^G:=\inf\left\{t\geq 0:\:\{X_0^G,\dots,X_t^G\}=V(G)\right\}
\end{equation}
to be the first time that $X^G$ has hit every vertex of $G$. We note that if
\begin{equation}\label{tildtdef}
\tilde{\tau}_{\rm cov}^G:=\inf\left\{t\geq 0:\:L^G_t(x)>0,\:\forall x\in V(G)\right\},
\end{equation}
then $\tilde{\tau}_{\rm cov}^G=\tau_{\rm cov}^G+1$; this equality will be crucial for our argument. We note that, as in the first part of this section, the steps we follow are not specific to the Sierpi\'nski gasket example, and will apply to any sequence of graphs for which we have a scaling limit for the random walks and local times. In order to state our main result, for a diffusion $X^F$ with state space $F$ and corresponding local times $(L_t^F(x))_{x\in F,t\geq0}$, we define $\tau_{\rm cov}^F$ and $\tilde{\tau}_{\rm cov}^F$ analogously to (\ref{tdef}) and (\ref{tildtdef}), respectively.

{\cor Let $(G_i)_{i\geq 0}$ be the sequence of Sierpi\'nski gasket graphs of Section \ref{nestedsec}. If $x_i\in V(G_i)$ satisfy $x_i\rightarrow x$, then
\begin{eqnarray}
\limsup_{i\rightarrow\infty}\mathbf{P}^{G_i}_{x_i}\left(5^{-i}\tau_{\rm cov}^{G_i}\leq t\right)&\leq &
\mathbf{P}^{F}_{x}\left(\tau_{\rm cov}^{F}\leq t\right),\label{lower}\\
\liminf_{i\rightarrow\infty}\mathbf{P}^{G_i}_{x_i}\left(5^{-i}\tau_{\rm cov}^{G_i}\leq t\right)&\geq &
\mathbf{P}^{F}_{x}\left(\tilde{\tau}_{\rm cov}^{F}< t\right),\label{upper}
\end{eqnarray}
for every $t\geq 0$, where $\mathbf{P}_x^F$ is the law of the Brownian motion $X^F$ on the Sierpi\'nski gasket $F$ started from $x$.}
\begin{proof} Suppose that $t<{\tau}_{\rm cov}^F$. Then there exists an $x\in F$ such that $x$ is not contained in the set $\{X_s^F:\: 0\leq s\leq t\}$. By the continuity of $X^F$, it follows that there exists an $\varepsilon>0$ such that $B_E(x,\varepsilon)\cap \{X_s^F:\: 0\leq s\leq t\}=\emptyset$, where $B_E(x,\varepsilon)$ is the Euclidean ball of radius $\varepsilon$ centred at $x$. Now, applying the Skorohod representation theorem, it is possible to assume that we have a realisation of the relevant processes such that the convergence at (\ref{scaling}) occurs almost-surely. Since we are assuming convergence in the uniform topology, it follows that, for large $i$, $B_E(x,\varepsilon/2)\cap \{X_{5^is}^{G_i}:\: 0\leq s\leq t\}=\emptyset$, and so $5^it\leq  {\tau}_{\rm cov}^{G_i}$. Thus we conclude that $\liminf_{i\rightarrow\infty}5^{-i}{\tau}_{\rm cov}^{G_i}\geq {\tau}_{\rm cov}^F$, and the bound at (\ref{lower}) follows.

Suppose that $t>\tilde{\tau}_{\rm cov}^F$. As local times are increasing in $t$, it must be the case that $L_t^F(x)>0$ for every $x\in F$. Together with the continuity of the local times, this implies that there exists an $\varepsilon>0$ such that $L_t^F(x)>\varepsilon$ for every $x\in F$. Again applying the Skorohod representation theorem, we may suppose that the conclusion of Theorem \ref{localscaling} holds almost-surely.
Since this statement is also in the uniform topology, it follows that $6(3/5)^iL_{5^it}^{G_i}(x)>\varepsilon/2$ for every $x\in V(G_i)$ for large $i$. It thus holds that $5^it\geq \tilde{\tau}_{\rm cov}^{G_i}$ for large $i$, which establishes $\limsup_{i\rightarrow\infty}5^{-i}{\tau}_{\rm cov}^{G_i}=\limsup_{i\rightarrow\infty}5^{-i}\tilde{\tau}_{\rm cov}^{G_i}\leq\tilde{\tau}_{\rm cov}^F$. This readily yields the statement at (\ref{upper}).
\end{proof}

{\rem \label{f2} (a) One can check that $0<{\tau}^F_{\rm cov}\leq\tilde{\tau}^F_{\rm cov}<\infty$, almost-surely (cf.~the proof of \cite[Theorem 6.3]{BP}), and so the limiting expressions are non-trivial.\\
(b) It is an interesting open problem to determine for which limiting diffusions the identity ${\tau}^F_{\rm cov}=\tilde{\tau}^F_{\rm cov}$ holds almost-surely, as it does for reflected Brownian motion on an interval, for example. Indeed, if this were true for the Brownian motion on the Sierpi\'nski gasket, then the above result would actually demonstrate that $5^{-i}{\tau}_{\rm cov}^{G_i}\rightarrow {\tau}^F_{\rm cov}$ in distribution.\\
(c) In fact, in the Sierpi\'nski gasket case, it is possible to check that $5^{-i}{\tau}_{\rm cov}^{G_i}\rightarrow {\tau}^F_{\rm cov}$ in distribution using the time-change argument of Remark \ref{f1}. Indeed, if
\[\tau_{\rm cov}^{F,i}:=\inf\left\{t\geq0:\:V(G_i)\subseteq\{X^F_s:\:0\leq s\leq t\}\right\},\]
then it is possible to check from the continuity of $X^F$ that $\tau_{\rm cov}^{F,i}\rightarrow\tau_{\rm cov}^{F}$, almost-surely. Hence we also have that $\tau_i(\tau_{\rm cov}^{F,i})\rightarrow \tau_{\rm cov}^{F}$, almost surely. Since $\tau_i(\tau_{\rm cov}^{F,i})$ is the cover time of the continuous time random walk with exponential mean $5^{-i}$ holding times, by a reparametrisation and the law of large numbers, it follows that the rescaled cover times of the discrete time walks converge in distribution. We reiterate, though, that we expect the argument of this article to be more widely applicable than this.}

\section*{Acknowledgements}

The author would like to thank Julia Komjathy for her comments on a draft version of this article that led to several improvements.

\bibliography{david}
\bibliographystyle{amsplain}

\end{document}